\documentclass[a4paper,10pt]{article}
\def\uu{\bigsqcup}
\def\eps{\varepsilon}

\usepackage[T2A]{fontenc}
\usepackage[cp1251]{inputenc}
\usepackage[english,russian]{babel}
\usepackage[tbtags]{amsmath}
\usepackage{amsfonts,amssymb}

\usepackage{mathrsfs}

\usepackage{graphicx}

\begin{document}
\Large

\title{ Slow convergences of the ergodic averages 
}
\author{Valery V. Ryzhikov}
\date{11.11.22}
\maketitle

\begin{abstract}
Birkhoff's theorem states that for an ergodic automorphism, the time averages converge to the space average.  Given sequence $\psi(n)\to+0$, U. Krengel proved that for  any ergodic automorphism there is  an  indicator  such  that   the    corresponding time averages  converged a.e.  with a rate  slower than $\psi$.  We  prove  again similar statements  answering  a question of I. Podvigin in passing.

\end{abstract}

\maketitle


\section{Введение }
Пусть $T$ -- эргодический автоморфизм  вероятностного пространства $(X,\mu)$
(автоморфизмы мы также называем преобразованиями). Для функции $f\in L_1(X,\mu)$ 
ее эргодические средние сходятся относительно  $L_1$-нормы  к константе $\int f$:
$$f_N:=\frac 1 N \sum_1^N T^if \to\int f \ \ (N\to\infty).$$ 
В статье будем использовать обозначение $f_N(T,x)$ вместо $f_N(x)$.
Теорема Биркгофа гласит, что для   п.в. $x\in X$ 
$$f_N(T,x):=\frac 1 N \sum_1^N T^if(x) \to\ \int_X f(x)d\mu.$$

Кренгелем \cite{K} был обнаружен интересный  эффект отсутствия оценок для скорости сходимости в эргодических теоремах. Мы дадим    доказательство результата Кренгеля о сходимостях по норме (Теорема 1) на основе теоремы Альперна (Теорема 2)  и ответим  на вопрос Подвигина из \cite{P} (Теорема 3).    Для множества произвольной меры из интервала $(0,1)$ и последовательности   $\psi(N)\to +0$ будет  построен эргодический автоморфизм $T$ такой, что  временные  средние для индикатора этого множества почти всюду сходятся медленнее $\psi(N)$.  Теорема Кренгеля \cite{K} 
универсальна в том смысле, что  для всякого эргодического автоморфизма и $\psi(N)\to +0$ она гарантирует  существование соответствующего множества меры $\frac 1 2$, для которого имеет место указанная медленная сходимость временных средних.

\vspace{2mm}
\bf Теорема 1.   \it Для  $\psi(N)\to +0$ и  эргодического автоморфизма $T$  найдется множество  $A$ такое, что для постоянной функции 
$a=\mu(A)$ выполняется 
$$\limsup_{N\to\infty}  \frac {\left\| f_N -a\right\|}{\psi(N)}=\infty.$$ \rm 

\vspace{2mm}
Ниже будет дано  доказательство, использующее обобщение леммы Рохлина-Халмоша, предложенное Альперном \cite{A} (см. также \cite{PR}). Башней высоты $h$ для автоморфизма $T$ называется набор непересекающихся множеств вида $T^iB$, $0\leq i<h$, или их объединение  $\uu_{i=0}^{h-1}T^iB$, $\mu(B)>0$.

\vspace{2mm}
\bf Теорема 2. \rm (Alpern) \it Пусть натуральные числа $h(j)$, $j=1,2,\dots$, взаимно просты в совокупности и $\sum_ja_j=1$,  $ a_j>0$.  Для апериодического $T$ найдется разбиение фазового пространства $X=\uu X_j$, где $X_j$ --  башня высоты $h(j)$, причем  $\mu(X_j)=a_j$. \rm

\vspace{2mm}
Доказательство теоремы 1.  Положим для опрделенности $a_j=2^{-j}$.  
Выбираем натуральные числа $h(j)$ так, чтобы для всех $j$ выполнялись условия 
$$\frac {2^{-j-1}}  {\psi(h(j))} \to\infty, \ \ \frac {\psi(h(j))} {\eps_jh(j)}\to\infty,$$
где 
$$\eps_j=\sum_{i>j}\frac {a_i}{h(i)}.$$
(Последовательность $h(j)$ можно выбирать так, чтобы $\frac {a_i}{h(i)}>2\eps_j$.)

 Обозначим через   $a$ постоянную функцию  $a(x)=\mu(X_1)=0.5$,
положим  $f=\chi_{X_1}$ и
$$U_j=\uu_{i\leq j} X_i, \ \ C_j=\uu_{i>j} X_i.$$

Заметим, что функция $f_{h(j)}$ на $C_j$ равна 0, 
кроме множества меры не большей,
чем  $\eps_jh(j)$. Последняя величина есть мера объединения первых $h_j$ этажей  в башнях $X_i$  с номерами $i>j$.
Из сказанного вытекает, что 
$$ \|f_{h(j)}\, \chi_{C_j}\| \leq \eps_jh(j),$$
$$  \int_{C_j}\,|\,f_{h(j)}-a\,|d\mu >0.5 \mu(C_j)-\eps_jh(j)> 2^{-j-2}>> \psi(h(j)),$$
$$ \frac {\|f_{h(j)}-a  \|} {\psi(h(j))} \to\infty.$$
Теорема доказана.

Следующее утверждение  отвечает на один из вопросов работы  \cite{P}.

\vspace{2mm}
\bf Теорема 3. \it Для всякой  $\psi(N)\to +0$ и числа $a\in (0,1)$ найдeтся множество $A$, $\mu(A)=a$,  и  эргодический автоморфизм $T$     такие, что для  индикатора $f=\chi_A$ и почти всех $x$ выполнено 
$$\limsup_{N\to\infty}  \frac {\left| f_N(T,x) -\mu(A)\right|}{\psi(N)}=\infty.$$ \rm 

\vspace{2mm}
Доказательство.
Пусть $X=\uu_j A_j$,  $j=0,1,2,\dots$,  $\mu(A_j)=a_j$, $\sum_ja_j=1$,  $ a_j>0$. 
 Положим
$$c = \mu(A_0), \ \ c_j=\frac {c}{1-\sum_{i>j}\mu(A_i)},\ \ j=1,2, \dots $$
($c_j$ --  среднее значение функции $f=\chi_{A_0}$ на множестве $A_0\cup A_1\cup \dots\cup A_j $).

Пусть задано преобразование $T_0$ такое, что  $T_0A_j=A_j$ для всех  $j$, причем ограничения $T_0$  на множества $A_j$ являются эргодическими.  Приступим к поиску
$T$. Фиксируем  последовательность $\eps_j\to +0$ такую, что  $\eps_j$  мало по сравнению с числом $(c_j-c)$.

\bf Шаг 1.  \rm Рассмотрим автоморфизм $T_1$, который  отличется от $T_0$ на множестве $\Delta_1$  маленькой меры,  $\Delta_1\subset A_0\cup A_1 $ и 
автоморфизм $T_1$ эргодичен 
на объединении $U_1=A_0\cup A_1 $. Существование такого $T_1$ для сколь угодно малой меры множества $\Delta_1$  очевидным образом вытекает, например, из теоремы 2. Действительно, изменяя апериодическое  преобразование (нас сейчас интересует только апериодическая часть  на $U_1$)  на самых последних этажах  башен, мы можем получить изоморфную копию любого другого апериодического преобразования, в частности, копии всех эргодических на $U_1$.  Но нам будет достаточно какого-нибудь  $T_1$. Отметим, что малость меры множества, на котором мы меняем преобразование, обеспечивается большой высотой башен. 

Так как автоморфизм $T_1$ эргодичен на $U_1=A_0\cup A_1 $, для $f=\chi_{A_0}$
в силу теоремы Биркгофа временные средние $f_{N}(T_1, x)$ сходятся к пространственному среднему $c_1$ п.в на $U_1$, а эта сходимость  влечет за собой сходимость по мере.   Таким образом,  для любого $\eps_1>0$ для всех достаточно больших $N_1$ выполнено 
 $$\mu\left(x\in U_1\,:\,  |f_{N_1}(T_1,x) -c_1|<\eps_1\right)>\mu(U_1)-
 {\eps_1}.$$

Зафиксируем такое число $N_1$ при условии, что $\psi(N_1)< c_1-c$ (сейчас это неравенство не играет принципиальной роли  -- мы просто  готовимся к описанию дальнейших шагов).

\bf Шаг 2.  \rm  Рассмотрим автоморфизм $T_2$, который  отличется от $T_1$ на множестве $\Delta_2$ (очень) маленькой меры,  $\Delta_2\subset A_1\cup A_2 $, и автоморфизм $T_1$ эргодичен на объединении $U_2=A_0\cup A_1 \cup A_2$.  (Изменяя
$T_1$, мы можем  обеспечить дизъюнктность  множеств  $\Delta_2$ не
и $\Delta_1$.)  Множество $\Delta_2$   мы выбираем столь малой меры, что 
при замене $T_1$ на $T_2$  по-прежнему выполняется 
$$\mu\left(x\in U_1\,:\,  |f_{N_1}(T_2,x) -c_1|<\eps_1\right)>\mu(U_1)-
{\eps_1}.$$

Из теоремы Биркгофа имеем:  для любого $\eps_2>0$ для всех достаточно больших $N_2$
 $$\mu\left(x\in U_2\,:\,  |f_{N_2}(T_2,x) -c_2|<\eps_2\right)>\mu(U_2)-{\eps_2}.$$
Зафиксируем такое число $N_2$ при условии, что 
$$\psi(N_2)< \frac {c_2-c} 2 \ .$$

\bf Шаг $\bf j+1$. \rm Рассмотрим автоморфизм $T_{j+1}$, который  отличется 
от $T_j$ на множестве $\Delta_{j+1}$. Выполнены условия:   $\Delta_{j+1}\subset A_j\cup A_{j+1} $, $\Delta_{j+1}$ не
пересекается с $\Delta_i$ при $i\leq j$,  а автоморфизм $T_{j+1}$ эргодичен 
на объединении $U{j+1}=A_0\cup A_1\cup\dots \cup A_{j+1}$. 
Множество $\Delta_{j+1}$    выбирается настолько малым, что для всех  $i\leq j$
выполняется 
$$\mu\left(x\in U_i\,:\,  |f_{N_i}(T_{j+1},x) -c_i|<\eps_i\right)>\mu(U_i)-
\eps_i.$$

В силу  теоремы Биркгофа   для любого $\eps_{j+1}>0$ и для всех достаточно больших $N_{j+1}$ имеем
 $$\mu\left(x\in U_{j+1}\,:\,  |f_{N_{j+1}}(T_{j+1},x) -c_{j+1}|<\eps_{j+1}\right)>\mu(U_{j+1})-
 {\eps_{j+1}}.$$
Зафиксируем такое число $N_{j+1}$ при условии, что 
$$\psi(N_{j+1})< \frac {c_{j+1}-c} {j+1}\,.\eqno (1)$$

Автоморфизмы $T_j$ очевидным образом сходятся к автоморфизму $T$,
для которого сохраняются все неравенства с той лишь оговоркой, что строгие неравенства формально надо заменить на  нестрогие: для всех $i$ имеем
$$\mu\left(\, x\in U_i\,:\,  |f_{N_i}(T,x) -c_i|\leq\eps_i\,\right)\geq\mu(U_i)-
\eps_i.\eqno (2)$$

Принимая во внимание $\mu(U_j)\to 1$, $\eps_j\to 0$, с учетом $(1),(2)$  мы получаем, что автоморфизм $T$ и  множество $A=A_0$ ($c=\mu(A)$) удовлетворяют утверждению 
$$\lim_{j \to\infty}  \frac {\left| f_{N_j}(T,x) -\mu(A)\right|}{\psi(N_j)}=\infty.$$
Теорема доказана.

\vspace{3mm}

\rm
\large

\end{document}